\newcommand{\be}{\begin{eqnarray}}
\newcommand{\ee}{\end{eqnarray}}
\newcommand{\re}[1]{(\ref{#1})}

\def\mb{\mbox{}}

\def\ll#1{\left#1}
\def\r#1{\right#1}

\def\p{\partial}

\def\bd{\begin{displaymath}}
\def\ed{\end{displaymath}}
\def\ba#1{\begin{array}{#1}}
\def\ea{\end{array}}
\def\nn{\nonumber}
\newfont{\symb}{msam7 scaled 1200}
\def\znakr{\raise1.5pt\hbox{\symb\char66\kern-2pt\char74}}
\def\znakl{\raise1.5pt\hbox{\symb\char73\kern-2pt\char67}}
\newfont{\Bbb}{msbm10 scaled 1200}
\documentstyle[12pt]{article}

\title{Projective representations of k--Galilei group}
\author{C.Gonera, P.Kosi\'nski\thanks{supported by KBN grant 2 P03B 130 12}, P.Ma\'slanka$^*$ \\
Department of Theoretical Physics II\\
University of \L\'od\'z\\
Pomorska 149/153\\
90--236 \L\'od\'z, POLAND\\
M.Tarlini\\
Dipartimento di Fisica, Universita di Firenze\\
and\\
INFN --- Firenze, Italy}
\date{}
\begin{document}
\maketitle
\begin{abstract}
The projective representations of k--Galilei group G$_k$\  are found by contracting
the relevant representations of $\kappa$--Poincare group. The projective multiplier
is found. It is shown that it is not possible to replace the projective representations
of G$_k$\ by vector representations of some its extension.
\end{abstract}
\section{Introduction}
There has been some attention paid recently to the deformations of space--time
symmetries depending on dimensionful parameter, the so--called $\kappa$--symmetries \cite{b1}--\cite{b11}.
They are interesting because they provide rather mild deformation of classical
space--time symmetries with dimensionful parameter (cut--off?) naturally built in.
It is, of course, still open question whether quantum symmetries provide a proper
way of introducing a fundamental energy/length scale into the theory;
in particular, special attention should be paid to the problems related
to noncocommutativity of coproduct which apparently seems to be in some
contradiction with kinematical properties of many--particle systems
(see, however, Ref.\cite{b9}). In spite of that it could be interesting
to study in more detail the properties of $\kappa$--deformed space--time
symmetries. Some preliminary studies of their physical implications
were already undertaken. In particular, Bacry \cite{b10} has found that they posses
some attractive features from the point of view of general requirements imposed
on kinematical symmetries.

In most papers that appeared so far the deformations of Poincare symmetry were
studied. However, it seems to be interesting to analyse the deformation
of its nonrelativistic counterpart, i.e.\ the deformed Galilei group. One version
of deformed Galilei algebra was studied in Ref.\cite{b4} where it was shown to 
provide the symmetry algebra of one--dimensional Heisenberg ferromagnet. In Ref.\cite{b3}
another Galilei algebra was found by applying the contraction procedure
($c\to\infty,\;\kappa\to 0,\;k\equiv\kappa c$\ kept fixed) 
to $\kappa$--Poincare
algebra in trigonometric version. The properties of this algebra (in hiperbolic
version) as well as the algebra obtained by letting 
$c\to\infty,\;\kappa\to \infty,\;k\equiv\kappa/c$\ -- fixed, were studied in Ref.\cite{b8}.
Finally, in Ref.\cite{b11}, the $\kappa$--Poincare group was contracted to the 
k--Galilei group and the latter was shown to be dual to k--Galilei algebra.
The bicross--product structure of both was revealed and the projective
representations of two--dimensional counterpart of k--Galilei group were constructed.

In the present paper we continue the study of k--Galilei group. In Sec.II the 
projective multiplier is found by contracting the trivial multiplier on $\kappa$--Poincare
group. In the $k\to\infty$\ limit it reduces to the standard nontrivial multiplier
on classical Galilei group. In Sec.III the unitary projective representations
 of k--Galilei group are constructed (again by contraction from the representations
of $\kappa$--Poincare group) and their infinitesimal form is given. The generators
of infinitesimal representations form the algebra which is a ``central" extension
of k--Galilei algebra. The question arises whether this structure can be lifted
to the Hopf algebra structure. This question is equivalent to the one whether
there exists a ``central" extension of k--Galilei group such that the projective
representations of the latter are equivalent to the standard (vector) representations
of its central extension. In Sec.IV we prove that no such central extension exists.
Sec.V is devoted to some conclusions. Finally, the technicalities are relegated
to the Appendices.

We conclude the introduction with short resume of results obtained in Ref.\cite{b11}.

In order to find the k--Galilei group one can apply the contraction procedure to
the $\kappa$--Poincare group defined in Ref.\cite{b2}. The following convenient
parametrization of Lorentz group can be used for contraction procedure (actually, it
differs slightly from the one adopted in Ref.\cite{b2})
\be\label{w1}
{\Lambda^0}_0&=&\frac{1}{\sqrt{1-\vec{v}^2/{c^2}}}\equiv \gamma\nn\\
{\Lambda^0}_i&=&\frac{\gamma}{c}v^i\nn\\
{\Lambda^i}_0&=&\frac{\gamma}{c}v^i\\
{\Lambda^i}_j&=&(\delta^i_k+(\gamma-1)\frac{v^iv^k}{\vec{v}^2}){R^k}_j\nn\\
RR^T&=&I\nn
\ee
as well as
\be
a^0&=&c\tau\nn
\ee
The following Hopf algebra G$_k$\ (k--Galilei group) is obtained from the contraction\ 
$\kappa\to 0,\;c\to\infty,\;k\equiv\kappa c$--fixed:
\be\label{w2}
\mb[{R^i}_j,{R^k}_l]=0,\;\;\;[{R^i}_j,v^k]=0,\;\;\;&[v^i,v^j]=0\nn\\
\mb[a^i,a^j]=0,\;\;\;[\tau,a^i]=\frac{i}{k}a^i\nn\\
\mb[\tau,v^i]=\frac{i}{k}v^i,\;\;\;[\tau,{R^i}_j]=0\nn\\
\mb[v^i,a^j]=\frac{i}{k}(\frac{1}{2}\vec{v}^2\delta_{ij}-v^iv^j)\nn\\
\mb[{R^i}_j,a^k]=\frac{i}{k}(\delta_{ik}v^m{R^m}_j-v^i{R^k}_j)\\
\Delta{R^i}_j={R^i}_k\otimes {R^k}_j\nn\\
\Delta v^i={R^i}_j\otimes v^j+v^i\otimes I\nn\\
\Delta a^i={R^i}_j\otimes a^j+v^i\otimes\tau+a^i\otimes I\nn\\
\Delta \tau=\tau\otimes I+I\otimes \tau&\nn\\
({R^i}_j)^*={R^i}_j,\;\;\;(v^i)^*=v^i,&(a^i)^*=a^i,\;\;\;\tau^*=\tau;\nn
\ee
G$_k$\ has a bicross--product structure
\be
G_k&=&T^*\znakr C(E(3))\nn
\ee
where C(E(3)) is the algebra of functions on classical group E(3) generated
by ${R^i}_j$\ and $v^i$\ while $T^*$\ is defined by
\be
\mb[\tau,a^i]=\frac{i}{k}a^i,\;\;\;[a^i,a^j]=0\nn\\
\Delta a^i=a^i\otimes I+I\otimes a^i,\;\;\;\Delta\tau=\tau\otimes I+I\otimes\tau\nn
\ee
The k-Galilei algebra $\widetilde{G_k}$, dual to G$_k$, reads 
\be\label{w3}
\mb[J_i,J_k]=i\epsilon_{ikl}J_l,\;\;\;[J_i,L_k]=i\epsilon_{ikj}L_j,\;\;\;[J_i,P_k]=i\epsilon_{ikj}P_j\nn\\
\mb[L_i,H]=iP_i,\;\;\;[L_i,P_j]=\frac{i}{2k}\delta_{ij}\vec{P}^2-\frac{i}{k}P_iP_j\nn\\
\Delta J_i=J_i\otimes I+I\otimes J_i\nn\\
\Delta H=H\otimes I+I\otimes H\\
\Delta L_i=I\otimes L_i+L_i\otimes e^{-\frac{H}{k}}-\frac{i}{k}\epsilon_{ijk}J_i\otimes P_k\nn\\
\Delta P_i=I\otimes P_i+P_i\otimes e^{-\frac{H}{k}}\nn\\
P_i^*=P_i,\;\;\;H^*=H,\;\;\;L_i^*=L_i,\;\;\;J_i^*=J_i\nn
\ee
$\widetilde{G_k}$\ has also bicrossproduct structure
\be
\widetilde{G_k}=T\znakl U(J,L)\nn
\ee
where U(J,L) is universal covering of Lie algebra e(3) while T is defined by
\be
\mb[H,P_i]=0,\;\;\;[P_i,P_j]=0\nn\\
\Delta H=H\otimes I+I\otimes H,\;\;\;\Delta P_i=P_i\otimes e^{-\frac{H}{k}}+I\otimes P_i\nn
\ee
The duality rules are the same as in classical case.
\section{Projective multiplier on G$_k$}
In analogy with the classical case one can define projective representation of 
a quantum group A acting in a Hilbert space H as a map $\rho: H\to H\otimes A$\ 
satisfying
\be\label{w4}
(\rho\otimes I)\circ\rho(\psi)=(I\otimes \omega)((I\otimes\Delta)\circ\rho(\psi))
\ee
where $\omega$\ is a unitary element of $A\otimes A$\ (projective multiplier) 
obeying suitable consistency condition \cite{b11}.

Two projective representations $\rho$\ and $\rho'$\ are called equivalent if 
there exists a unitary element $\zeta\in A$\ such that
\be\label{w5}
\tilde{\rho}=(I\otimes\zeta)\rho
\ee
The corresponding multipliers are related by the formula
\be\label{w6}
(\zeta\otimes\zeta)\omega=\omega'\Delta(\zeta)
\ee
Obviously, a multiplier $\omega$\ is trivial (the representation is equivalent to the
vector one) if
\be\label{w7}
\omega=(\zeta^{-1}\otimes\zeta^{-1})\Delta(\zeta)
\ee
In the classical case it is sometimes possible to obtain nontrivial multiplier
by contraction \cite{b12}. Assume the group $\widetilde{G}$\ is obtained
from G by  contraction. Even if G does not admit nontrivial projective
multipliers, one can proceed as follows. Let $\zeta(g)$\ be a unitary 
function on G, $\zeta(g)\zeta^*(g)=1$. Define a trivial multiplier on G
\be\label{w8}
\omega(g,g')=\zeta^*(g)\zeta^*(g')\zeta(gg')
\ee
It can happen that, $\zeta(g)$\ being properly chosen, the specific combination
of $\zeta$'s appearing on the right hand  side of eq.\re{w8}\ tends to the well defined limit under
contraction while $\zeta(g)$\ itself has no such a limit. We can then expect 
that the limiting $\omega(g,g')$\ is a nontrivial multiplier on $\widetilde{G}$.
This is, for example, the case for G being Poincare group and
\be\label{w9}
\zeta(\{\Lambda,a\})=e^{-imca^0}
\ee
The corresponding multiplier $\omega(g,g')$, eq.\re{w8}, gives in the contraction limit $c\to\infty$
\be\label{w10}
\tilde{\omega}=e^{-im(\frac{\vec{v}^2}{2}\tau'+v^k{R^k}_ia'^{i})}\;\;\;,
\ee
the standard multiplier on Galilei group.

Following the classical case we define the trivial projective multiplier on $\kappa$--Poincare
group:
\be\label{w11}
\omega=(e^{imca^0}\otimes e^{imca^0})e^{{-imc}({\Lambda^0}_\mu\otimes a^\mu+a^0\otimes I)}
\ee

Our aim is to find the limiting form of $\omega$\ for $\kappa\to 0,\;\;c\to\infty,\;\;k\equiv \kappa c$;
the $\kappa,c$--dependence of m is yet unknown and must be determined from condition
that the nontrivial limit exists. To this end we rewrite first $\omega$\ in a
more convenient form making explicit the cancellation of divergent terms. The long
and rather tedious calculations reported in Appendix lead to the following expression
for $\omega$
\be\label{w12}
\omega=e^{i(mc-\kappa\ln(ch(\frac{mc}{\kappa})+{\Lambda^0}_0sh(\frac{mc}{\kappa})))\otimes a^0}e^{-i\kappa\frac{sh(\frac{mc}{\kappa}){\Lambda^0}_k}{ch(\frac{mc}{\kappa})+{\Lambda^0}_0sh(\frac{mc}{\kappa})}\otimes a^k}
\ee 
In order to calculate the limiting value of $\omega$\ we use the parametrization \re{w1}. It is easy
to check that in order to obtain the nontrivial limit one can choose the following form of m:
\be\label{w13}
m=\frac{-k}{2c^2}\ln(1-\frac{2Mc^2}{k})
\ee
where M is some fixed mass parameter (which we assume to be positive); let us note
that for k negative, $m>0$\ and $m\to 0$\ for $c\to\infty$\ while $mc^2\tau$\ diverges. 
Taking the $c\to\infty$\ 
limit in eq.\re{w12} one obtains
\be\label{w14a}
\tilde{\omega}=e^{-ik\ln(1+\frac{M\vec{v}^2}{2k})\otimes \tau}e^{-\frac{Mv^k{R^k}_i}{1+\frac{M\vec{v}^2}{2k}}\otimes a^i}
\ee
which can be also written as
\be\label{w14b}
\tilde{\omega}=e^{-i(\frac{2k}{M\vec{v}^2}\ln(1+\frac{M\vec{v}^2}{2k})\otimes I)(\frac{M\vec{v}^2}{2}\otimes\tau+Mv^k{R^k}_i\otimes a^i)}
\ee
This expression is a natural generalization of the one obtained in Ref.\cite{b11} for two--dimensional
case. In the classical limit $k\to\infty$\ it coincides with the standard multiplier
\re{w10}.

Let us conclude this section by noting one trouble related to the formula \re{w14a}. In
order to keep the Poincare mass m real we had to assume k negative. This, however,
implies that $\tilde{\omega}$\ is singular somewhere. On the other hand, with k positive,
$\tilde{\omega}$\  is everywhere regular.

It seems that this trouble cannot be cured in a simple way. The following argument
can be given to support this point of view. Using the results contained in Ref.\cite{b18}
one easily concludes that the general form of irreducible (co)representation of
the $\kappa$--Poincare group is obtained by replacing the exponentials on the 
right--hand side of eq.\re{w15} by
\be
e^{-i\kappa\ln(\frac{p_0+C}{A})\otimes a^0}e^{-\frac{i\kappa p_k}{p_0+C}\otimes a^k}\nn
\ee
where $A=A(m,c,\kappa),\;C=C(m,c,\kappa)$\ are two real functions subject to the
condition $C^2-A^2=m^2c^2$\ but otherwise arbitrary. So the question arises whether
our trouble can be cured by an appropriate choice of $A$\ and $C$\ such that obtains,
in the $c\to\infty,\;\kappa\to 0$\ limit, the representation given by formula~(19)
while $m=m(M,c,k)$\ lies, for $k>0$, in the physical region $m>0$. This  seems not
to be possible. Let us put again $\vec{p}=\frac{m}{M}\vec{q}$\ and consider the second
exponential. It is easy to see that, in order to obtain the proper limiting formula,
the following condition should be fulfiled
\be
\lim\limits_{c\to\infty}c^2(1+\frac{C}{mc})=\frac{k}{M}\nn
\ee
However, due to the condition $C^2-A^2=m^2c^2,\;|C|\geq mc$\ and the above equation
can be satisfied only provided $C=-mc-\Delta,\;\Delta\geq 0$. Then
\be
\lim\limits_{c\to\infty}\frac{c\Delta}{m}=-\frac{k}{M}\nn
\ee
which is impossible for $k>0,\;M>0$\ and $m>0$.

\section{Contraction of representations}
The unitary representations of the $\kappa$--Poincare group were constructed in
Ref.\cite{b13} (see also \cite{b14}). This constructions can be summarized as follows.
The representation space is the Hilbert space of square integrable (with respect to
the standard measure $d^3\vec{p}/2p_0$) functions over the hiperboloid $p^2=m^2$\ 
taking their values in the vector space carrying the spin s representation of 
rotation group (s is assumed to be integer, for s halfinteger one should consider
quantum ISL(2,\Bbb C\rm)\ group \cite{b15} which only amounts to small modifications).
 The (right) corepresentation reads
\be\label{w15}
\rho:f_i(p_\mu)\to e^{-i\kappa\ln(ch\frac{mc}{\kappa})+\frac{p_0}{mc}sh(\frac{mc}{\kappa}))\otimes a^0}e^{\frac{-i\kappa sh(\frac{mc}{\kappa})p_k}{mc ch(\frac{mc}{\kappa})+p_0sh(\frac{mc}{\kappa})}\otimes a^k}\cdot\nn\\
\mb\cdot D_{ij}(R(p\otimes I,I\otimes\Lambda))f_j(p_\nu\otimes{\Lambda^\nu}_\mu);
\ee
here by $D(R(p\otimes I,I\otimes \Lambda))$\ we denote the spin s representation of standard
Wigner rotation written as an element of tensor product $H\otimes A$.

It follows from eq.\re{w15} that the whole deformation is contained in translation
sector; in other words the representation is obtained by integrating the 
infinitesimal representation given in Majid--Ruegg basis \cite{b16},\cite{b17}.
In the limit $\kappa\to\infty$\ unitary representations of classical Poincare
group are recovered.

In order to find the representations of k--Galilei group we apply again the
contraction procedure. To this end we put  $\kappa=k/c$\ and take 
$m\equiv m(M,k,c)$\ as defined by eq.\re{w13}. As in the classical case it is
necessary to subtract the rest energy by redefining $\rho$:
\be\label{w16}
\tilde{\rho}\equiv (I\otimes e^{-imca^0})\rho
\ee

Finally, contrary to the  classical case, we have to redefine the momenta and
the wave functions as follows:
\be\label{w17}
\vec{p}=\frac{m}{M}\vec{q},\;\;\;\frac{m}{M}f_i(\vec{p})=\tilde{f_i}(\vec{q})
\ee

It is now easy to check that the limit $c\to\infty$\ exists and gives the following
unitary representation of k--Galilei group:
\be\label{w18}
\tilde{\rho_{nr}}:\tilde{f_i}(\vec{q})\to e^{-ik\ln(1+\frac{\vec{q}^2}{2Mk})\otimes \tau}e^{\frac{-iq_k}{1+\vec{q}^2/2Mk}\otimes a^k}(I\otimes D_{ij}(R))\tilde{f_j}(q_k\otimes{R^k}_i+\nn\\
\mb+I\otimes Mv^k{R^k}_i)
\ee
acting in the Hilbert space of functions square integrable with respect to invariant
measure $d^3\vec{q}$\ and taking values in the vector space carrying spin s
representation of rotation group.

As a next step let us find the infinitesimal form of representation $\rho_{nr}$.
Let us recall that if 
\be
\rho: H\ni f\to f_{(\alpha)}\otimes a_{(\alpha)}\in H\otimes A \nn
\ee
is the (right) corepresentation of the quantum group A then any element X of
the dual Hopf algebra (quantum Lie algebra) is represented by the operator
\be\label{w19}
\tilde{X}: H\ni f\to f_{(\alpha)}<a_{(\alpha)},X>\in H
\ee

The relevant duality rules can be, as it was mentioned above, adopted from classical
theory. A simple calculation then gives
\be\label{w20}
J_k=-i\epsilon_{klm}q_l\frac{\p}{\p q_m}+s_k\nn\\
L_k=iM\frac{\p}{\p q_k}\nn\\
H=k\ln(1+\frac{\vec{q}^2}{2Mk})\\
P_k=\frac{q_k}{1+\frac{\vec{q}^2}{2Mk}}\nn
\ee
Let us note that H and P$_k$\ are nonsingular only provided $k>0$.

The operators \re{w20} verify the following commutation rules
\be\label{w21}
\mb[J_i,J_k]=i\epsilon_{ikl}J_l,\;\;\;[J_i,P_k]=i\epsilon_{ikl}P_l,\;\;\;[J_i,L_k]=i\epsilon_{ikl}L_l\nn\\
\mb[K_i,H]=iP_i\\
\mb[K_i,P_j]=iM\delta_{ij}e^{-\frac{2H}{k}}+\frac{i}{2k}\delta_{ij}\vec{P}^2-\frac{i}{k}P_iP_j\nn
\ee
For $M\to 0$\ this algebra coincides with the algebraic sector of k--Galilei algebra~\re{w3}.

Finally, let us note the following dispersion relation valid within the representation~\re{w20}
\be\label{w22}
k(1-e^{-\frac{H}{k}})=\frac{\vec{P}^2}{2M}
\ee
\section{The central extension of k--Galilei group}
It is well known that in the classical case given a projective representation
of a group G one can construct the group G' such that this projective representation
of G is equivalent  to the vector
representation of G'. The natural question arises whether the analogous construction
is possible in quantum case.

Let us assume that there exists a Hopf algebra $G'_k$\ such that
\begin{description}
\item[(i)] $G_k'$\ is obtained from $G_k$\ by adding one new unitary element
$\zeta: \zeta\zeta^*=\zeta^*\zeta=I$
\item[(ii)] $G_k$\ is a Hopf subalgebra of $G_k'$
\item[(iii)] $\Delta(\zeta)=(\zeta\otimes\zeta)\omega$
\end{description}
where $\omega$\ is a projective multiplier.

Then, if $\rho$\ is a projective representation of $G_k$\ determined by $\omega$\ 
(cf.\ eq.\re{w4}),
\be
\rho'=(I\otimes\zeta)\rho\nn
\ee
is a vector representation of $G_k'$.

In the commutative case (iii) determines $G_k'$\ uniquely and consistently.
In the quantum case, however, $\Delta$\ should be a homomorphism which, together
with (iii), imposes nontrivial consistency conditions. It has been already
shown \cite{b8} that $G_k'$\ cannot be obtained by straightforward generalization
of standard contraction from trivial extension of $\kappa$--Poincare. We show below
that there is no solution to the problem, at least if the existence of well--defined
limit $k\to\infty$\ which reproduces the classical situation is assumed. To this
end let us note first that eqs.\re{w2} and \re{w14a} imply
\be\label{w23}
e^{\alpha(\tau\otimes I+I\otimes\tau)}\tilde{\omega}e^{-\alpha(\tau\otimes I+I\otimes\tau)}=e^{i\ll(\frac{2k}{M\vec{v}^2}\ln(1+\frac{M\vec{v}^2}{2k}e^{\frac{2i\alpha}{k}})\otimes I\r)\ll(\frac{M\vec{v}^2}{2}\otimes\tau+Mv^k{R^k}_i\otimes a^i\r)}
\ee
Therefore
\be\label{w24}
\mb[\tau\otimes I+I\otimes\tau,\tilde{\omega}]=\frac{2}{k}\tilde{\omega}(M\vec{v}^2\otimes\tau+Mv^k{R^k}_i\otimes a^i)(\frac{1}{1+\frac{M\vec{v}^2}{2k}}\otimes I)
\ee
The homomorphism condition
\be\label{w25}
\mb[\Delta(\zeta),\Delta(\tau)]=\Delta([\zeta,\tau])
\ee
gives
\be\label{w26}
\Delta(\zeta)\ll(\frac{2}{k}(\frac{1}{1+\frac{M\vec{v}^2}{2k}}\otimes I)(\frac{M\vec{v}^2}{2}\otimes\tau+Mv^k{R^k}_i\otimes a^i)\r)+\nn\\
\mb+(I\otimes\zeta)([\zeta,\tau]\otimes I)\tilde{\omega}+(\zeta\otimes I)(I\otimes [\zeta,\tau])\tilde{\omega}&=&\Delta([\zeta,\tau])
\ee
It follows from eq.\re{w26} that the commutator $[\zeta,\tau]$\ should be of the form
\be\label{w27}
\mb[\zeta.\tau]=\frac{2}{k}\zeta X
\ee
where X is an element of G$_k$, the 1/k factor is extracted out explicitly and the factor 2 is written
for convenience.

The following relation follows immediately from eqs.\re{w26} and \re{w27}
\be\label{w28}
\Delta(X)=\tilde{\omega}^{-1}(X\otimes I)\tilde{\omega}+\tilde{\omega}^{-1}(I\otimes X)\tilde{\omega}+(\frac{1}{1+\frac{M\vec{v}^2}{2k}}\otimes I)\cdot\nn\\
\mb\cdot(\frac{M\vec{v}^2}{2}\otimes\tau+Mv^k{R^k}_i\otimes a^i)
\ee
Taking the lowest term in 1/k expansion we get
\be\label{w29}
\frac{M\vec{v}^2}{2}\otimes\tau+Mv^k{R^k}_i\otimes a^i=\Delta(X)-X\otimes I-I\otimes X
\ee
which can be viewed as the relation on classical Galilei group. However, eq.\re{w29} does not
hold true because the left--hand side is not a coboundary. 
\section{Conclusion}
Using the contraction technique we have found the projective multipliers on k--Galilei
group and the corresponding projective representations, both in global as well as 
in infinitesimal form. It appears that we obtain a well--defined and regular structure
for $c\to\infty$\ provided the deformation parameter k is taken to be positive.
On the other hand, in order to keep the Poincare mass parameter in the allowed 
region in the course of contraction we should rather assume k to be negative.
We do not have a clear understanding of this phenomenon.

In the classical case the projective representations can be always converted into the 
vector representations of suitably defined extension of the original group.
We have seen in Sec.IV that this is not necessarily the case for quantum
groups. There exists no suitable extension of $G_k$\ which, in the classical
limit $k\to\infty$\, reduces to the standard case. This seems to be not a serious
obstacle because it concerns some technical rather then fundamental aspect
of the theory.

The problem which certainly deserves further study is the multiplication of
representations. This is important if we would like to reconcile the 
noncommutativity of the algebra coproduct with the basic properties of 
many--particle systems, especially those containing identical particles.

\appendix
\section*{Appendix}
\setcounter{equation}{0}
\renewcommand{\theequation}{A.\arabic{equation}}
We derive here eq.\re{w12}. In order to simplify the notation we omit the tensor
product symbol $\otimes$\ writing instead with prime the factors appearing right
to it. Eq.\re{w11} reads then
\be\label{wa1}
\omega=e^{imca^0}e^{imca'^0}e^{-imc({\Lambda^0}_0a'^0+{\Lambda^0}_ka'^k+a^0)}
\ee
Define
\be\label{wa2}
X(m)\equiv e^{imca^0}e^{-imc({\Lambda^0}_0a'^0+{\Lambda^0}_ka'^k+a^0)}
\ee
Then X(0)=I and X(m) obeys the equation
\be\label{wa3}
\stackrel{\bullet}{X}(m)=(-ic)(Y_0(m)a'^0+Y_k(m)a'^k)X(m)
\ee
where
\be\label{wa4}
Y_\mu(m)=e^{imca^0}{\Lambda^0}_\mu e^{-imca^0},\;\;\;Y_\mu(0)={\Lambda^0}_\mu\;\;\;
\ee
Let us first calculate $Y_0(m)$. Using $\kappa$--Poincare group commutation rules
we get
\be\label{wa5}
\stackrel{\bullet}{Y}_0(m)=ice^{imca^0}[a^0,{\Lambda^0}_0]e^{-imca^0}=-\frac{c}{\kappa}e^{imca^0}(({\Lambda^0}_0)^2-1)e^{-imca^0}=\nn\\
\mb=-\frac{c}{\kappa}(Y_0^2(m)-1)
\ee
As $Y_0(m)$\ commute for all m, \re{wa5} together with the initial condition \re{wa4} can
be solved by separation of variables yielding
\be\label{wa6}
Y_0(m)=\frac{{\Lambda^0}_0ch(\frac{mc}{\kappa})+sh(\frac{mc}{\kappa})}{{\Lambda^0}_0sh(\frac{mc}{\kappa})+ch(\frac{mc}{\kappa})}
\ee
With $Y_0(m)$\ explicitly known one can apply the same procedure to find $Y_k(m)$;
the result reads
\be\label{wa7}
Y_k(m)=\frac{{\Lambda^0}_k}{ch(\frac{mc}{\kappa})+{\Lambda^0}_0sh(\frac{mc}{\kappa})}
\ee
Now we can go back to the eq.\re{wa3}. It cannot be solved directly because two
terms on the right--hand side do not commute. To account for this we pass to ``interaction
picture" and define
\be\label{wa8}
X(m)=
e^{-ic\int\limits^m_0dmY_0(m)a'^0}W(m),\;\;\;W(0)=I
\ee
For \re{wa3} we get
\be\label{wa9}
\stackrel{\bullet}{W}(m)=-icY_k(m)e^{ic\int\limits^m_0dmY_0(m)a'^0}
a'^ke^{-ic\int\limits^m_0dmY_0(m)a'^0}W(m)
\ee
On the other hand
\be\label{wa10}
ic\int\limits^m_0dmY_0(m)=i\kappa\ln(ch(\frac{mc}{\kappa})+{\Lambda^0}_0sh(\frac{mc}{\kappa}))
\ee
It follows from \re{wa9} and \re{wa10} that
\be\label{wa11}
\stackrel{\bullet}{W}(m)=\frac{-ic{\Lambda^0}_ka'^k}{(ch(\frac{mc}{\kappa})+{\Lambda^0}_0sh(\frac{mc}{\kappa}))^2}W
\ee
Again everything commutes so that \re{wa11} can solved
\be\label{wa12}
W=e^{-ic\int\limits^m_0\frac{dm}{(ch(\frac{mc}{\kappa})+{\Lambda^0}_0sh(\frac{mc}{\kappa}))^2}{\Lambda^0}_ka'^k}
\ee
or
\be\label{wa13}
W=e^{i\kappa\frac{sh(\frac{mc}{\kappa}){\Lambda^0}_ka'^k}{ch(\frac{mc}{\kappa})+{\Lambda^0}_0sh(\frac{mc}{\kappa})}}
\ee
Eq.\re{w11} follows directly from \re{wa1}, \re{wa2}, \re{wa8} and \re{wa13}. 

\end{document}